\documentclass[leqno,twoside]{article}
\usepackage[centertags]{amsmath}
\usepackage{theorem,euscript,eufrak,amscd,verbatim}
\newcommand{\sh}[1]{\EuScript{#1}}     
\newcommand{\gothic}[1]{\EuFrak{#1}}   
\newcommand{\mb}[1]{\mathbf{#1}}       
\newcommand{\fn}[1]{\mathrm{#1}}       
\newcommand{\ctop}{c_{\mathrm{top}}}   
\newcommand{\Co}{\mb{C}}               
\newcommand{\ext}{\sh{E}xt}            
\newcommand{\hiso}{\ensuremath{\sim}}  
\newcommand{\isoa}{\xra{\sim}}          
\renewcommand{\O}{\sh{O}}              
\newcommand{\red}{_{\text{red.}}}
\renewcommand{\hom}{\sh{H}om}          
\newcommand{\T}{\sh{T}}                
\newcommand{\Z}{\mb{Z}}                
\newcommand{\viso}{\textstyle {\wr}}   
\newcommand{\abs}[1]{\ensuremath{\lvert #1\rvert}}
\newcommand{\C}[2]{C_{#1\backslash #2}}           
\newcommand{\Hilb}[1][{}]{\ensuremath{\mathrm{Hilb}^{#1}}}   
\newcommand{\I}[2]{\sh{I}_{#1\backslash #2}}      
\newcommand{\N}[2]{\sh{N}_{#1\backslash #2}}      
\renewcommand{\P}[1][{}]{\ensuremath{\mb{P}^{#1}}}
\newcommand{\rest}[2]{\left. #1\right|_{#2}}      
\newcommand{\shd}[1]{#1^\vee}                     
\DeclareMathOperator{\Bl}{Bl}
\DeclareMathOperator{\CH}{CH}
\DeclareMathOperator{\h}{h}
\renewcommand{\H}{\ensuremath{\operatorname{H}}}
\DeclareMathOperator{\Ext}{Ext}
\DeclareMathOperator{\Hom}{Hom}
\DeclareMathOperator{\Obs}{Obs}
\DeclareMathOperator{\Pic}{Pic}
\DeclareMathOperator{\rProj}{\mathbf{Proj}}  
\newcommand{\R}{\mathrm{R}}
\DeclareMathOperator{\rank}{rk}
\DeclareMathOperator{\Res}{Res}
\DeclareMathOperator{\Spec}{Spec}
\DeclareMathOperator{\rSpec}{\mathbf{Spec}}   
\DeclareMathOperator{\Sym}{Sym}
\newcommand{\hra}{\hookrightarrow}
\newcommand{\iso}{\cong}
\newcommand{\lra}{\longrightarrow}

\newcommand{\tensor}{\otimes}

\newcommand{\wt}[1]{\widetilde{#1}}
\newcommand{\xra}[1]{\xrightarrow{#1}}
\renewcommand{\epsilon}{\varepsilon}

\newcommand{\blank}{\underline{\hspace{1em}}}
\newcommand{\E}{\sh{E}}
\renewcommand{\i}{\sh{I}}
\renewcommand{\j}{\sh{J}}
\newcommand{\K}{\sh{K}}
\renewcommand{\L}{\sh{L}}
\newcommand{\m}{\gothic{m}}
\newcommand{\p}{\gothic{p}}

\newcommand{\Ccomp}{\gothic{C}}
\newcommand{\Q}{\sh{Q}}
\newcommand{\Sg}{\Xi}
\newcommand{\V}{\sh{V}}

\makeatletter
\renewcommand{\subsection}{%
\@startsection{subsection}{2}{0mm}{-\baselineskip}%
{-0.5em}{\normalfont\normalsize\bfseries}%
}
\renewcommand{\subsubsection}{%
\@startsection{subsubsection}{3}{0mm}{-0.5\baselineskip}%
{-0.5em}{\normalfont\normalsize\scshape}%
}
\renewcommand{\@seccntformat}[1]{\csname the#1\endcsname\hspace{0.5em}}
\makeatother

\newcommand{\subsub}[1]{%
\setcounter{subsubsection}{\value{equation}}%
\subsubsection{#1}%
\setcounter{equation}{\value{subsubsection}}%
}

\renewcommand{\thesubsubsection}%
  {(\arabic{section}.\arabic{subsection}.\arabic{subsubsection})} 

\theoremstyle{change}
\theoremheaderfont{\normalfont\bfseries}
{\theorembodyfont{\rmfamily\slshape}
\newtheorem{theorem}{Theorem}[section]
\newtheorem{proposition}[theorem]{Proposition}
\newtheorem{lemma}[theorem]{Lemma}
\newtheorem{corollary}[theorem]{Corollary}}

{\theorembodyfont{\normalfont\normalshape}}

\newenvironment{thm}%
{\setcounter{theorem}{\value{subsection}}\begin{theorem}}%
{\end{theorem}\setcounter{subsection}{\value{theorem}}}

\newenvironment{lem}%
{\setcounter{theorem}{\value{subsection}}\begin{lemma}}%
{\end{lemma}\setcounter{subsection}{\value{theorem}}}

\newenvironment{cor}%
{\setcounter{theorem}{\value{subsection}}\begin{corollary}}%
{\end{corollary}\setcounter{subsection}{\value{theorem}}}

\newenvironment{prop}%
{\setcounter{theorem}{\value{subsection}}\begin{proposition}}%
{\end{proposition}\setcounter{subsection}{\value{theorem}}}

\newcommand{\qedsymbol}{$\clubsuit$}
\newcommand{\qed}{\hspace{\fill}\qedsymbol}
\numberwithin{equation}{subsection}

\title{Counting curves which move with threefolds}
\author{Herbert Clemens and Holger P. Kley}
\date{July 1998}

\begin{document}
\addtocounter{section}{-1}
\maketitle

\begin{abstract}
Let $X$ be a (possibly nodal) $K$-trivial threefold moving in a fixed
ambient space $P$.  Suppose $X$ contains a continuous family of curves, all
of whose members satisfy certain unobstructedness conditions in $P$.  A
formula is given for computing the corresponding virtual number of
curves, that is, the number of curves on a generic 
deformation of $X$ ``contributed by'' the continuous family on~$X$.
\end{abstract}

\section{Introduction}

\subsection{} Suppose $X_0$ is a projective threefold with at worst
ordinary node singularities which is embedded in a smooth
projective manifold~$P$, and that
\[
Z_0\subseteq X_0
\]
is a connected curve.  Assume that
\[
\H^1(Z,\N{Z}{P})=0
\]
for all $\{Z\}\in J'$, an open set in the Hilbert scheme of~$P$, with
$J'$ containing the connected component $I'$ of $\{Z_0\}$ in the Hilbert
scheme of~$X_0$.  Suppose in addition
that
\[
\omega_{X_0}\tensor\O_Z \iso \O_Z
\]
for all $\{Z\}\in I'$.  Then the expected dimension of the set of
curves $\{Z\}\in I'$ which deform to a generic deformation of $X_0$ is
zero.  The purpose of this paper is to compute the (virtual) number
$\gamma(I')$ of such curves under certain additional assumptions.

We work in the setting in which $X_0$ is the zero-scheme of a regular
section of a vector bundle on $P$ such that, via pull-back and
push-forward, $I'$ is given as the zero-scheme of the associated
section $\sigma_0$ of the associated bundle $V$ on~$J'$.  Then
\[\gamma(I') :=\deg(\ctop(V)),\]
which can be computed as 
the geometric intersection number of $\sigma_0(J')$ with the
zero-section $z_{V}(J')$ of $V$.  Using Fulton-MacPherson intersection
theory~\cite{Fulton_IT}, one rescales $\sigma_0(J')$ by multiplication
by larger and larger constants to create a ``homotopy'' between
$\sigma_0(J')$ and the normal cone $\C{I'}{J'}\subset \rest{V}{I'}$ of
$I'$ in $J'$, so that 
\[
\gamma(I') = z_V(J')\cdot\C{I'}{J'}
\]
can be calculated as an intersection product in~$\rest{V}{U}$.

Under the above assumptions, we reinterpret the sheaf of obstructions to
deformation as the sheaf of K\"ahler differentials
on~$I'$ with logarithmic poles along the locus of curves passing
through the nodes of~$X_0$, thereby allowing the computation
of~$\gamma(I')$ in terms of the geometry of~$I'$.

\medskip

The results of the present paper do not guarantee that $I'$
contributes $\gamma(I')$ rigid curves to a general deformation $X_t$
of~$X_0$;  {\sl a priori}, the count $\gamma(I')$ is purely virtual.
In some cases, however, the precise structure of the obstruction sheaf
does enable rigidity results;  see \cite{Kley98a} and the remarks
following example~\ref{ex:ciCY} of the present work.

Furthermore, our computations are \emph{local} in that we work with a
single connected component of the Hilbert scheme of ~$X_0$.  The
global (virtual) number of curves of given degree and genus is 
computed---at least in case $g=0,1$---via Gromov-Witten invariants in
the fundamental papers of Kontsevich \cite{Kontsevich95}, Givental
\cite{Givental97, Givental98} and Lian, Liu, Yau~\cite{LLY97}.  Also,
a symplectic treatment of Gromov-Witten invariants for nodal $X_0$
appears in~\cite{LR98}.

\medskip

The paper is organized as follows:  In \S\ref{sec:generalities} we
establish some general properties of Hilbert schemes, allowing $X_0$
to have arbitrary isolated singularities and $Z$ and $X_0$ to be of
arbitrary positive dimension.  Let $d=\dim Z$; then by Serre duality,
$\H^d(\shd{\N{Z}{X_0}}\tensor\omega_Z)$ is isomorphic to
$\shd{\H^0(\N{Z}{X_0})}$ which is the cotangent space to $I'$ at~$\{Z\}$.
The key result---Lemma~\ref{lem:1.9}---establishes a relative form of
this fact:  $\Omega^1_{I'}$ is the $d$th
higher direct image of the relative conormal tensor the relative
dualizing sheaf.  The assumptions and notation established
in \S\ref{sec:generalities} will be used throughout the work.  

Then in \S\ref{sec:smooth}, we explore the
case in which $X_0$ is a smooth threefold, $Z\subset X_0$ a curve, and
$\omega_{X_0}\tensor\O_Z \iso \O_X$.  These assumptions give an
adjunction isomorphism between $\H^1(\N{Z}{X_0})$ and
$\H^1(\shd{\N{Z}{X_0}}\tensor\omega_Z)$, 
which in light of the results of \S\ref{sec:generalities}, gives an
isomorphism between $\Omega^1_{I'}$ and the first higher direct image
of the relative normal bundle; this is
Proposition~\ref{prop:OmegaI'}.  Given the role of $\H^1(\N{Z}{X_0})$
in the  obstruction theory of the Hilbert scheme, the resulting
formula for $\gamma(I')$ in Corollary~\ref{cor:Xsmooth} should not be
too surprising.  The main technical difficulty arises from the failure
of $I'$ to be smooth in general.

In \S\ref{sec:nodes}, we extend the computation of \S\ref{sec:smooth}
to the case in 
which $X_0$ has ordinary nodes, which---at least for enumerative
purposes---introduces logarithmic poles along the locus of curves
passing through the nodes of~$X_0$.  We will assume that the generic
curve parameterized by $I'$ does not pass through the nodes of~$X_0$
which enables us to use the technical results of the previous section.

Finally, in \S\ref{sec:exas}, we illustrate our formulas with three
examples.  

\subsection{Conventions and Notation}
All schemes are separated and of
finite type over the field $\Co$ of complex numbers.
If $Z$ is a closed subscheme of $X$, we denote by $\I{Z}{X}$ the ideal
sheaf, by $\C{Z}{X}:= \rSpec(\bigoplus
\I{Z}{X}^n\big/\I{Z}{X}^{n+1})\to Z$ the {\sl normal cone} of $Z$
in~$X$, and by $\N{Z}{X} := \hom_{\O_X}(\I{Z}{X},\O_Z) =
\hom_{\O_Z}(\I{Z}{X}\big/\I{Z}{X}^2,\O_Z)$ the {\sl normal sheaf} of 
$Z$ in~$X$.  If $X$ is smooth, $\T_X:=\shd{\left(\Omega^1_X\right)}$
is the {\sl tangent sheaf} of~$X$.  When appropriate, the notations
$\omega_X$ or $\omega_{X/X'}$ refer to dualizing or relative dualizing
sheaves, not canonical sheaves.

We generally pass between the notions of locally-free sheaf and vector
bundle without comment, but we use roman type $E\to X$ to denote the
geometric vector bundle associated to a locally-free
$\O_X$-module~$\E$.  In this case, $z_E:X\to E$ is the zero-section
of~$E$.

For flat families of schemes we use the notation $F\to F'$;  that is,
a prime (${}'$) indicates the base of a flat family.

Finally if $Z\hra X$ is a closed embedding, we denote by $\{Z\}$ the
corresponding point in the Hilbert scheme of~$X$.

\subsection{Acknowledgments}  We wish to express our gratitude to
Y.~Ruan for explaining to us his symplectic methods of computing
Gromov-Witten invariants of Calabi-Yau threefolds with nodes, and to
L.~Ein for suggesting the connection with the 
log-complex.  Those discussions and the the previous results of the
second author for elliptic curves on Calabi-Yau complete intersections
were the genesis for the results in this paper.  We also thank the
referee for several corrections and suggestions for improving the
exposition. 

\section{Generalities}\label{sec:generalities}

\subsection{}\label{ss:P}Let $P$ be a projective manifold, let $d\ge
1$, and let $J'$ be a 
connected open subscheme of the Hilbert scheme of $P$, maximal with
respect to the properties that:

\subsub{}\label{ass:lci} For all $\{Z\}\in J'$, $Z$ is a connected
local complete intersection scheme of dimension $d$ in $P$.

\subsub{}\label{ass:H1Nvanishing} For all $\{Z\}\in J'$, the normal
bundle $\N{Z}{P}$
satisfies
\[
\H^1(Z,\N{Z}{P}) = 0.
\]
So $J'$ is a smooth, irreducible quasi-projective variety.

\subsection{}\label{ss:X}Suppose further that we are given a locally
free sheaf $\E$ on $P$ and a regular section
\[
s_0\in\H^0(P,\E)
\]
such that the zero scheme
\begin{equation}
X_0:=(s_0=0)\subset P
\end{equation}
of $s_0$ has isolated singularities.  We require that $\E$ be
sufficiently ample, in the sense that for all $\{Z\}\in J'$,
\begin{equation}\label{eq:H1Evanishing}
\H^1(Z,\E\tensor\O_Z)=0.
\end{equation}

Let $E\to P$ denote the geometric vector bundle associated to $\E$.
There is a natural surjection
\[
\I{z_E(P)}{E} \to \I{X_0}{s_0(P)}\iso\I{X_0}{P}
\]
which restricts to give a sequence of isomorphisms:
\begin{equation}\label{eq:Edual}
\begin{split}
\shd{\E}\tensor\O_{X_0} &\iso
\left(\I{z(P)}{E}\big/\I{z(P)}{E}^2\right)\tensor\O_{X_0}\\
 &\iso \I{X_0}{s_0(P)}\big/\I{X_0}{s_0(P)}^2\\
 &\iso \I{X_0}{P}\big/\I{X_0}{P}^2.
\end{split}
\end{equation}

\subsection{}  Next consider the incidence scheme or universal
family $J\subset J'\times P$ with projections
\begin{equation*}
\begin{CD}
         J @>{q}>> P \\
      @VV{p}V               @.\\
         J'         @.      {}
\end{CD}
\end{equation*}
Let
\begin{equation*}
\V := p_*q^*\E \quad\text{and}\quad\sigma_0:=p_*q^*s_0.
\end{equation*}
Then because of \eqref{eq:H1Evanishing} (see
\cite[Thm.\@~1.5]{Kley98}), $\V$ is locally free and there is a
scheme-theoretic equality
\begin{equation}\label{eq:zero_scheme}
I':= \Hilb[X_0]\cap J' = (\sigma_0=0).
\end{equation}
Notice that that by shrinking $J'$, we may assume $I'$ is connected.

\subsection{}\label{ss:functoriality}   For the remainder of this
section, consider an 
arbitrary cartesian square
\begin{equation}\label{eq:B}
\begin{CD} B @>{b}>> J \\
@VV{p^B}V @VV{p}V\\
B' @>{b'}>> J'
\end{CD}
\end{equation}
where we allow $B'$ to be an object in the analytic category.  Setting 
\[
A':= B'\times_{J'} I',
\]
\eqref{eq:B}  pulls back to a cartesian square
\[
\begin{CD} A @>{a}>> I \\
@VV{p^A}V @VV{p^0}V\\
A' @>{a'}>> I'
\end{CD}
\]
Let
\[
q_0:=\rest{q}{I}\colon I \to X_0
\]
be the second projection, and set
\[
q_B :=q\circ b\quad\text{and}\quad q_A:=q_0\circ a.
\]

\subsection*{} In these situations, we often will need the ideal
sheaves
\[
\j_B:=\I{B}{B\times P}\quad\text{and}\quad \i_A:=\I{A}{A'\times X_0}.
\]
In the cases $B'=J'$ (and $A'=I'$), we simplify to
\begin{equation*}
\j:=\I{J}{J'\times P}\quad\text{and}\quad\i := \I{I}{I'\times X_0}.
\end{equation*}

\begin{lem}\label{lem:1.8}
There is an isomorphism of exact sequences of $\O_A$-modules:
\[\begin{CD}
0 @>>>q_B^*\shd{\E}\tensor\O_{A} @>{q_B^*s_0}>>
 b^*\left(\j\big/\j^2\right) \tensor\O_A @>>>
 a^*\left(\i\big/\i^2\right) @>>> 0\\
@. @| @VV{\viso}V @VV{\viso}V @.\\
0 @>>>q_B^*\shd{\E}\tensor\O_{A} @>{q_B^*s_0}>>
 \left(\j_B\big/\j_B^2\right) \tensor\O_A @>>>
 \left(\i_A\big/\i_A^2\right) @>>> 0\\
\end{CD}\]
\end{lem}

\subsection{Proof:}  Let $Z$ be a fiber of $A\to A'$, which we identify
with its image $Z\subset X_0$.  In light of 
\eqref{eq:Edual}, either left-hand arrow restricts on $Z$ to the standard
morphism of conormal sheaves.
\[
\I{X_0}{P}\big/\I{X_0}{P}^2\tensor\O_Z \to \I{Z}{P}\big/\I{Z}{P}^2.
\]
We assumed that $Z$ is a local complete intersection in the smooth
variety $P$, so away from the singularities of $X_0$, this is injective.
So, since $\dim Z\ge 1$ and $X_0$ has isolated singularities, it is
generically injective.  Suppose the kernel is supported at
some~$z\in Z$.

Since $\I{X_0}{P}\big/\I{X_0}{P}^2\tensor\O_Z$ is locally free, we may
restrict to an affine neighborhood $\Spec R$ of $z$ over which the kernel
corresponds to a non-zero submodule $M\subset R^{\oplus k}$ of a free
module and assume that $M$ is annihilated by the maximal ideal $\m$
of~$z$. Then any projection of $M$ to the various summands is
annihilated by~$\m$, and one of these must be non-zero, so is some
non-zero ideal $\pi(M)\subset R$ with proper support~$z$.  Therefore 
$\m$ is an (embedded) associated prime of~$R$.  But $Z$, being a
local complete intersection, is Cohen-Macaulay, and hence has no
embedded components.  (See, e.g., \cite[Thm.\@~17.3]{Matsumura_CRT}.)   
We conclude that the the morphism is everywhere injective.

Next observe that the natural map
\[
\I{z(P)}{E} \lra \I{X_0}{s_0(P)}\iso \I{X_0}{P}
\]
pulls back to a surjection
\[
q_0^*\I{z(P)}{E} \lra q_0^*\I{X_0}{P}
\]
on $I'\times P$, whence an exact sequence
\[
q_0^*\I{z(P)}{E}\lra \j  \lra \i \lra 0
\]
which remain exact after tensoring with $\O_Z$. 
Thus, both sequences restrict to exact sequences on each fiber of
$p^A$, so by Nakayama's lemma, they are exact.

Finally, the natural surjection $b^*\j\to \j_B$ induces a surjection
of the middle terms.  Since both are locally free of the same rank,
it must be an isomorphism.  It follows that the natural map between
the right-hand terms must also be an isomorphism.\qed

\subsection{} Let $\omega:=\omega_{A/A'}$ be the relative
dualizing sheaf.  The first two non-trivial sheaves in
Lemma~\ref{lem:1.8} are locally free $\O_A$-modules.  Therefore, by
Verdier duality, 
\begin{equation}\label{eq:Verdier1}
\begin{split}
(b')^*\shd{\V}\tensor\O_{A'} &=
\hom_{A'}\left(p^A_*q_A^*\E,\O_{A'}\right)\\ 
                 &\iso \R^d\left(p^A_* \circ%
                       \hom_A(\blank,\omega)\right)
                        \left(q_A^*\E\right)\\  
                 &=\R^dp^A_*\left(q_A^*\shd{\E} 
                       \tensor\omega\right) 
\end{split}
\end{equation}
and
\begin{equation}\label{eq:Verdier2}
\begin{split}
(b')^*\Omega^1_{J'}\tensor\O_{A'} &=%
                    \hom_{A'}((b')^*\T_{J'}\tensor\O_{A'},\O_{A'})\\ 
                 &\iso \R^d\!\left(p^A_* \circ%
                     \hom_A(\blank,\omega)\right)%
                     \left(\hom_A\left(b^*\left(\j\big/\j^2\right)\tensor\O_A,%
                      \O_A\right)\right)  \\
                 &=\R^d p^A_*\left(b^*\left(\j\big/\j^2\right) \tensor%
                     \omega\right),
\end{split}
\end{equation}
where we have used the infinitesimal properties of the Hilbert
scheme to make the identification
$p_*\left(\shd{\left(\j\big/\j\right)}\right)\iso \T_{J'}$.
We combine these calculations:
\begin{lem}\label{lem:1.9}  There is a commutative diagram
\[
\newcommand{\tws}{\tensor\omega}
\minCDarrowwidth1pc
\begin{CD}
\R^d p^A_*(q_B^*\shd{\E}\tws) @>{q_B^*s_0}>>
    \R^dp^A_*(b^*(\j\big/\j^2)\tws) @>>> 
    \R^d p^A_*(a^*(\i\big/\i^2)\tws) @>>> 0\\
@VV{\viso}V  @VV{\viso}V @VV{\viso}V @.\\
(b')^*\shd{\V}\tensor\O_{A'} @>{(b')^*\sigma_0}>>
(b')^*\Omega^1_{J'}\tensor \O_A @>>> (a')^*\Omega^1_{I'} @>>> 0
\end{CD}
\]
with exact rows and all vertical maps isomorphisms.
\end{lem}

\subsection{Proof:}
Exactness of the top row results from Lemma~\ref{lem:1.8} after
applying $\R^d p^A_*$ (which preserves right exactness because the
fiber-dimension of $p^A$ is~$d$).

Next note that since $I'$ is the zero-scheme of $\sigma_0$, there is a
surjection
\[
\shd{\V}\xra{\sigma_0}\I{I'}{J'}
\]
and hence, using the standard exact sequence of differentials, a
commutative diagram with exact rows and columns:
\begin{equation}\label{eq:1.9.1}
\begin{CD}
\shd{\V}\tensor\O_{I'} @>>> \Omega^1_{J'}\tensor\O_{I'} @>>>
  \Omega^1_{I'} @>>> 0\\
@VVV   @| @| @.\\
\I{I'}{J'}\big/\I{I'}{J'}^2 @>>> \Omega^1_{J'}\tensor\O_{I'} @>>>
  \Omega^1_{I'} @>>> 0\\
@VVV @. @. @. \\
0 @. {} @. {} @. {}
\end{CD}
\end{equation}
Now the  exactness of the bottom row of the Lemma follows by pull-back 
to~$A'$.

Next the left two vertical isomorphism are 
\eqref{eq:Verdier1} and \eqref{eq:Verdier2}, and commutativity
of the left-hand square follows via duality, Nakayama's lemma, and
pull-back to $A'$ from~\cite[Prop.\@~1.6]{Kley98}.
Finally, a diagram chase establishes the existence of the third
vertical isomorphism and the commutativity of the right-hand square.
\qed

\subsection{}For the derived functors associated to the functor
\begin{equation}
\fn{T^A}:=p^A_*\circ\hom_A(\blank,\O_A)
\end{equation}
we have by \ref{ass:H1Nvanishing} that
\[
\R^1\fn{T}^A\left(b^*\left(\j\big/\j^2\right)\tensor\O_A\right) = 0.
\]
So from Lemma~\ref{lem:1.8} we obtain the exact sequence
\begin{multline}\label{eq:RT}
0\lra\fn{T}^A\left(a^*\left(\i\big/\i^2\right)\right) \lra%
 \left(b'\right)^*\T_{J'}\tensor\O_{A'} \lra\\
 \left(b'\right)^*\V\tensor\O_{A'}\xra{\delta}%
 \R^1\fn{T}^A\left(a^*\left(\i\big/\i^2\right)\right)\lra 0.
\end{multline}
We will see that, when $d=1$, the sheaf $\R^1\fn{T}^I(\i/\i^2)$ measures
the obstruction to moving the curves $\{Z\}\in I'$ when the regular
section $s_0$ of $\E$ is deformed to a generic section~$s_t$.

\subsection{}Let
\[
V:=\rSpec\left(\Sym^*(\shd{\V})\right)
\]
be the geometric vector bundle associated to~$\V$.  Define the {\sl
normal cone}
to $I'$ in $J'$ to be
\begin{equation*}
\begin{split}
\C{I'}{J'} &:=
\rSpec\bigoplus_{r=0}^{\infty}\left(\I{I'}{J'}^r\big/\I{I'}{J'}^{r+1}\right)
\end{split}
\end{equation*}
There is a canonical surjection
\[
\left(b'\right)^*\left(\I{I'}{J'}\right)\to\I{A'}{B'},
\]
so from the left-hand square of \eqref{eq:1.9.1}, we obtain a
commutative diagram of $\O_{A'}$-algebras
\begin{equation}\label{eq:graded}
\begin{CD}
\Sym^*\left(\left(b'\right)^*\shd{\V}\tensor\O_{A'}\right) @>>>
  \Sym^*\left(\left(b'\right)^*\Omega^1_{J'}\tensor\O_{A'}\right)\\ 
                 @VVV                             @|\\
\left(b'\right)^*\left(\bigoplus_{r=0}^\infty\left(\I{I'}{J'}^r\big/\I
  {I'}{J'}^{r+1}\right)\right) @>>> 
  \Sym^*\left(\left(b'\right)^*\Omega^1_{J'}\tensor\O_{A'}\right)\\
@VVV @.\\
\bigoplus_{r=0}^\infty\left(\I{A'}{B'}^r\big/\I
  {A'}{B'}^{r+1}\right) @. {}
\end{CD},
\end{equation}
whence a commutative diagram of morphisms of cones over $A'$:
\begin{equation}\label{eq:cones}
\begin{CD}
{} @. \C{A'}{B'}\\
@. @VVV\\
\rest{\left(b'\right)^*\left(T_{J'}\right)}{A'} @>>> A'\times_{I'}\C{I'}{J'}\\
          @|                @VVV\\
\rest{\left(b'\right)^*\left(T_{J'}\right)}{A'} @>>>
 \rest{\left(b'\right)^*V}{A'} 
\end{CD}.
\end{equation}
(See \cite[Ch.~4]{Fulton_IT} and \cite[\S1]{BF97}.)  Note that the
left-hand vertical arrows in \eqref{eq:graded} are surjective, so that
the right-hand vertical arrows in \eqref{eq:cones} are closed embeddings.

\subsection{}  Now \eqref{eq:zero_scheme} is precisely the statement
that the square
\[\begin{CD}
I' @>>> J'\\
@VV{i}V @VV{\sigma_0}V\\
J' @>{z_V}>> V
\end{CD}\]
is Cartesian.  Then we define the class
\begin{equation}
\gamma(I') := \Z(\sigma_0) := z_V^![J']\in \CH_*(I')
\end{equation}
to be the {\sl localized top Chern class} of \cite[\S14.1]{Fulton_IT},
where it is shown that 
\[
i_*\gamma(I')=\ctop(V)\cap[J']\in \CH_*(J'),
\]
i.e., that this class represents the top Chern class of~$V$.  In the
language of \cite{BF97}, this class is the {\sl virtual fundamental
class} 
\[
\gamma(I') = [I', F^\bullet],
\]
where $F^\bullet$ is the complex $\rest{[\shd{V}}{I'}\to
\rest{\Omega^1_{J'}}{I'}]$ (in degrees $-1$ and $0$);  see `the basic
example' in[loc.\@ cit., \S6].  See  also~\cite[\S1]{GP97}.
We remark that although the definition of $\gamma(I')$ in no way
depends on the completeness of~$I'$, its enumerative significance
does.

\section{Smooth threefolds}\label{sec:smooth}
\subsection{}\label{ss:3fold}It is at this point that we make our
final assumptions, namely assume:
\begin{equation}\label{eq:curves}
d = \dim Z =1
\end{equation}
so that in particular, since $\E\tensor\O_{X_0} \iso \N{X_0}{P}$,
\ref{ass:H1Nvanishing} implies \eqref{eq:H1Evanishing}.
Furthermore, we assume that
\begin{equation}\label{eq:threefold}
\dim X_0 = 3
\end{equation}
and that for all $\{Z\}\in I'$,
\begin{equation}\label{eq:Ktrivial}
\omega_{X_0}\tensor\O_Z \iso \O_Z\quad\text{and}\quad\h^0(\O_Z) = 1.
\end{equation}
Finally, we assume that
\begin{equation}\label{eq:Xsmooth}
\text{$X_0$ is smooth},
\end{equation}
although this will be weakened in \S\ref{sec:nodes}.

\subsection{}From \eqref{eq:Ktrivial} we immediately have
\begin{equation}\label{eq:K}
p^0_*q_0^*(\omega_{X_0}) =: \K
\end{equation}
for some line bundle $\K$ on~$I'$,  and we set 
\[
\K_{A'}:=(a')^*\K \xra{\hiso}  p^A_*q_A^*(\omega_{X_0}).
\]

Using (\ref{eq:curves}--\ref{eq:Xsmooth}), \ref{ass:H1Nvanishing} and
Riemann-Roch, we compute:
\begin{equation}
\begin{split}
\dim J' &= \chi(\N{Z}{P})\\
        &= \chi(\N{Z}{X_0}) + \chi(\E\tensor\O_Z)\\
        &= \deg \omega_{Z} + 2(1 - g_a(Z)) + \rank V\\
        &= \rank V.
\end{split}
\end{equation}
Consequently,
\[\gamma(I')\in \CH_0(I').\]

\subsection{}\label{sec:adjunction}  Now \eqref{eq:curves},
\eqref{eq:threefold}, and \eqref{eq:Xsmooth} imply
that $(\i\big/\i^2)$ is locally free of rank two, 
so the relative adjunction isomorphism 
\[
\omega_{A/A'}\tensor a^*\Lambda^2(\i\big/\i^2) \iso q_A^*(\omega_{X_0})
\]
induces an isomorphism:
\begin{equation}\label{eq:rel_adj}
\begin{split} 
\omega_{A/A'}\tensor (\i_A\big/\i_A^2) 
& \iso \hom_A\left((\i_A\big/\i_A^2),%
\omega_{A/A'}\tensor\Lambda^2(\i_A\big/\i_A^2)\right)\\  
&\iso  \hom_A\left(\i_A\big/\i_A^2,q_A^*(\omega_{X_0})\right)\\
&\iso  \hom_A\left(\i_A\big/\i_A^2,(p^A)^*(\K_{A'})\right).
\end{split}
\end{equation}
Thus, by the projection formula and Lemma~\ref{lem:1.9}, 
\begin{prop}\label{prop:OmegaI'}  There is an isomorphism
\[
\R^1 p^A_*\left(\hom_A(\i_A\big/\i_A^2,\O_A)\right)\tensor\K_A \iso
(a')^*\Omega^1_{I'}.
\]
\end{prop}

\subsection{}A further consequence of the local freeness of $\i_A/\i_A^2$
is that we have an isomorphism of derived functors
\begin{equation}
\R p^A_*\left(\hom_A(\i_A\big/\i_A^2,\O_A)\right) \isoa
 \R\fn{T}^A(\i_A\big/\i_A^2),
\end{equation}
which, when combined with \eqref{eq:RT} and
Proposition~\ref{prop:OmegaI'}, yields the exact sequence
\begin{equation}\label{eq:delta}
(b')^*\T_{J'}\tensor\O_{A'} \lra (b')^*\V\tensor\O_{A'} \xra{\delta}
\hom_{A'}(\K_{A'},(a')^*\Omega^1_{I'}) \lra 0.
\end{equation}

\subsection{}\label{ss:Cnotation}  Let $\Ccomp$ be the set of
components of $\C{I'}{J'}$.  For $C\in\Ccomp$, let $m(C)$ be the
geometric multiplicity of $C$ in~$\C{I'}{J'}$, $S'=S'(C)\subset I'$
its support, $\p = \p(C)\in S'\subset I'$ the generic point
of~$S'$, $k= k(S')$ the function field of~$S'$, and
$C_\p:=C\times_{I'} \p$ the fiber of $C$ over~$\p$.   Following our
conventions, $p^S\colon S \to S'$ is the pullback to $S'$ of
$p^0\colon I \to I'$. 

Now the surjective sheaf morphisms 
\[
\V\tensor\O_{I'} \xra{\delta} \hom_{I'}(\K,\Omega^1_{I'})
  \xra{\mathrm{rest.}} \hom_{I'}(\K,\Omega^1_{S'}).
\]
give rise a surjective composition
\begin{equation*}
\V\tensor k \to \Hom_k(\K\tensor k, \Omega^1_{I'}\tensor k) \to
\Hom_k (\K\tensor k,\Omega^1_{S'}\tensor k)
\end{equation*}
of maps of $k$-vector spaces, which we can view 
as a morphism
\begin{equation}\label{eq:kvec}
\rest{V}{\p} \to \rest{(\shd{T}_{S'}\tensor K_{S'}^{-1})}{\p}
\end{equation}
of varieties over $k$
.

\begin{lem}\label{lem:Cp}
Assume that $p^S$ is generically smooth.  Then given $\kappa\in
\H^0(\omega_{X_0})$, the composition
\[
C_\p \xra{\text{\eqref{eq:cones}}} \rest{V}{\p}
\xra{\text{\eqref{eq:kvec}}} \rest{(\shd{T}_{S'}\tensor K_{S'}^{-1})}{\p}
\xra{(p_0)_*q_0^*(\kappa)} \rest{\shd{T}_{S'}}{\p} 
\]
is the constant map to zero.
Furthermore, if $p^0_*q_0^*(\kappa)$ does not vanish at~$\p$, then
$C_\p$ is the geometric kernel of the vector-space surjection
$\rest{V}{\p}\to\rest{(\shd{T}_{S'}\tensor K_{S'}^{-1})}{\p}$. 
\end{lem}

\subsection{Proof:}  The second assertion follows from the first since
\eqref{eq:kvec} is surjective and
\[
\dim_k C_\p + \dim S' = \dim J' = \rank V.
\]

To prove the the first assertion, let $\pi'\colon \wt{J}'\to J'$ be
the blow-up of $J'$ along~$I'$.  Then the epimorphism
\[
\bigoplus_{j=0}^\infty \I{I'}{J'}^j \lra \bigoplus_{j=1}^{\infty}
 \left(\I{I'}{J'}^j\big/\I{I'}{J'}^{j+1}\right) 
\]
of graded $\O_{J'}$-algebras induces---via $\rProj$--- a closed
embedding
\[
\P(\C{I'}{J'}) \hra \wt{J}'
\]
over $I'\hra J'$, with (scheme-theoretic) image the
exceptional divisor.

Now let $U\subset S'$ be a small analytic neighborhood of a general point in $S'$ and
let $R$ be any line bundle over $U$ which is a sub-cone
of~$\rest{C}{U}$;  in other words, $R$ is a {\sl ray}
in~$\rest{C}{U}$.  Now such an $R$~determines (and is determined by) a
section   
\[
\phi\colon U\to \rest{\P(C)}{U}
\]
of the projection $\rest{\P(C)}{U}\to U$.
Let $0\in\Delta$ be a one-dimensional disk with
parameter~$t$. Set 
\[
Y':=\Delta\times U.
\]
Shrinking $U$ as necessary, we can construct an embedding
\[
g'=g'_R\colon Y' \to \wt{J}'
\]
such that
\[
g'(0,u) = \phi(u)\quad\text{for all $u\in U$},
\]
and
\[
g'\left(\Delta^*\times U\right) \subset
J'\setminus I'\subset \wt{J}'.
\]

Let $f' := \pi'\circ g'\colon Y' \to J'$, and
consider the scheme
\[
W':= (b\circ f)^{-1}(I').
\]
Then $W'\iso \Spec\left(\Co[t]/(t^m)\right) \times U$, where~$m=m(C)$.
Observe that under the canonical identification $U\iso \{0\}\times
U = W'\red$, we have isomorphisms
\begin{equation}\label{eq:reduced}
R \iso (\C{W'}{Y'})\red \iso \C{W'}{Y'}\times_{W'} W'\red.
\end{equation}

Keeping the notation of \ref{ss:functoriality}, we have a fiber square
\[
\begin{CD}
W:=W'\times_{I'} I @>{w}>> I\\
@VVV @VV{p}V\\
W' @>{w'}>> I'
\end{CD}
\]
and we set $\i_W := \I{W}{W'\times X_0}$.
Then applying the functoriality in Lemma~\ref{lem:1.8}---once with $B'
= Y'$, once with $B'=W'\red$---we see that the exact sequence
\begin{equation*}
(f')^*\T_{J'}\tensor\O_{W'} \lra (f')^*\V\tensor\O_{W'} \xra{\delta}
\hom_{W'}(\K_{W'},(w')^*\Omega^1_{I'}) \lra 0.
\end{equation*}
of \eqref{eq:delta} restricts on $W'\red\iso U$ to
\begin{equation*}
\begin{CD}
\T_{J'}\tensor\O_{U} @>>> \V\tensor\O_{U} @>{\rest{\delta}{U}}>>
\hom_{U}(\K_{U},\Omega^1_{I'}\tensor \O_U) @>>> 0\\
@. @. @VV{\viso}V @.\\
{} @. {} @. \R^1p^U_*(\N{U\times_I' I}{U\times X_0}) @. {}
\end{CD}.
\end{equation*}

Let $\sh{R}$ be the sheaf of sections of~$R$.  Then starting with
\eqref{eq:reduced}, we have a sequence of morphisms
\begin{multline*}
\sh{R}\iso\C{W'}{Y'}\times_{W'} U \xra{\text{\eqref{eq:cones}}}
 \V\tensor\O_U \xra{\rest{\delta}{U}} \R^1 p^U_*(\N{U\times_{I'}
 I}{U\times X_0}) \\
\xra{\mathrm{rest.}}\hom_U(\K_U,\Omega^1_U) \xra{p_*q^* \kappa}
 \Omega^1_U. 
\end{multline*}
Then the image of $\Gamma(\sh{R})$ in $\Gamma\left(\R^1
p^U_*\left(\N{U\times_{I'} I}{U\times X_0}\right)\right)$ consists of the
obstructions to extending the family $W\to W'$ to a family of
subschemes of $X_0$ over the base $\Spec\Co[t]/(t^{m+1}) \times U$, so
that in the  language of \cite[Ch.\@~2]{Clemens98}, this image contains
only sections {\sl associated to the normal cone}.  Essentially, these
are the obstructions to curvilinear deformation.  Shrinking $U$ as 
needed,  Cor.~2.6 of [loc.\@cit.] with $r=d=1$ states that the
composition 
\[ 
\sh{R} \to \V\tensor\O_U \to \Omega^1_U
\]
is zero.  (It is here that the assumption on the generic smoothness of
$p^S$ becomes necessary.)
But by construction, the map $R\to \rest{V}{U}$ factors through $C$,
and since $U\subset S'$ is general and $R$ is an arbitrary ray in~$C$, 
the Lemma follows.\qed

\begin{cor}\label{cor:Xsmooth}
Suppose that for all $C\in\Ccomp$, the support $S'(C)$ is smooth and that
the generic curve $I_{\p(C)}$ is smooth.  Suppose further that, for
each $S'=S'(C)$ there is a section
\[
\kappa_S\in\H^0(\omega_{X_0})
\]
such that $(p_0)_*q_0^*(\kappa_S)$ does not vanish identically
on~$S'$.  Then
\[
\gamma(I') = \sum_{C\in\Ccomp} m(C) \ctop\left(\K_{S'}^{-1}\tensor
 \Omega^1_{S'(C)}\right)\cap[S'(C)].
\]
\end{cor}

\subsection{Proof:}  By the excess intersection formula and the
linearity of the intersection product, 
\[
\gamma(I') = \sum_{C\in \Ccomp} m(C) \deg (z_V^! [C]).
\]
Now since each $S'$ is smooth, Lemma~\ref{lem:Cp} and the surjectivity
of \eqref{eq:kvec} imply, by dimension,
that there is an exact sequence of vector bundles over~$S'$
\[
0 \lra C \lra \rest{V}{S'} \lra \shd{T}_{S'}\tensor K_{S'}^{-1} \lra 0
\]
so that
\[
z_V^![C] = \ctop(\Omega^1_{S'}\tensor \K_{S'}^{-1})\cap[S'].
\]
The desired formula follows immediately.\qed

\section{Threefolds with nodes}\label{sec:nodes}

\subsection{} In this section, we make the same assumptions as in
in \ref{ss:P}, \ref{ss:X} and \ref{ss:3fold}, except that we weaken
\eqref{eq:Xsmooth} and assume instead that
\begin{equation}
\text{the singularities of $X_0$ are a set $\Sg$ of ordinary double points.} 
\end{equation}
Moreover, we require that
\begin{equation}\label{eq:Obs}
\Obs(Z,X_0) \subseteq
  \H^1\left(Z,\shd{\left(\I{Z}{X_0}\big/\I{Z}{X_0}^2\right)}\right) 
\end{equation}
for every curve $Z=q_0(I_{y'})$ with $y'\in I'$, where
\[
\Obs(Z,X_0)\subseteq \Ext^1_Z\left(\I{Z}{X_0}\big/\I{Z}{X_0}^2,\O_Z\right)
\]
is the space of obstructions generated by (obstructed) curvilinear
deformations of $Z$ in~$X_0$.  (See \cite[p.\@ 29ff.]{Kollar_RC}.)
Note that away from points $y'\in I'$ which represent curves passing
through nodes of~$X_0$, the techniques of \S\ref{sec:smooth} apply.
In particular, away from such points, the reduction to the use of the
curvilinear obstruction space in the following proceeds just as in the
proof of Lemma~\ref{lem:Cp}. 

\subsection{}  We continue to use the notation established in
\ref{ss:functoriality} and \ref{ss:Cnotation}.  For each $C\in\Ccomp$, assume

\subsub{}\label{ass:Ssmooth} The support $S'=S'(C)$ is smooth and the
morphism 
\[
p^S\colon S:=S'\times_{I'} I \to S'\] 
is generically smooth.

\subsub{}\label{ass:kappa} There is a section
\[
\kappa_S\in\H^0(X_0,\omega_{X_0})
\]
such that $p^S_*q_S^*(\kappa_S)$ does not vanish identically on~$S'$.

\subsub{}\label{ass:Y} If
\[
q_S(S)\cap \Sg \ne \emptyset,
\]
then $X_0$ contains a surface $Y_S$, smooth at 
\[
Y_S\cap \Sg=:\{x^1_S,\dots,x^{r(S)}_S\}
\]
such that
\[
q_S(S)\subseteq Y_S.
\]

\subsub{}\label{ass:D} The scheme
\[
D_{S'}^i = p^S(q_S^{-1}(x^i_S))
\]
is smooth divisor in $S'$ for all $i$.

\subsub{}\label{ass:nc}
Either\begin{enumerate}
\item
$D_{S'}:= \sum_{i=1}^r D_{S'}^i$ is a normal-crossing divisor on~$S'$,
or
\item $I'$ is smooth connected, so that $\Ccomp=\{C\}$ with $S'(C)=
I'$.
\end{enumerate}

\begin{thm}\label{thm:nodes} Under assumptions
\ref{ass:Ssmooth}--\ref{ass:nc} 
\[
\gamma(I') =
 \sum_{C\in\Ccomp}m(C)\left(\ctop\left(\Q_{S'(C)}\right)\cap[S'(C)]\right), 
\]
where each $\Q_{S'}$ is a locally free sheaf which is an extension
\[
0 \lra \K_{S'}^{-1}\tensor\Omega^1_{S'} \lra \Q_{S'} \lra
 \bigoplus_{i=1}^{r(S)} \O_{D_{S'}^i}\lra 0,
\]

\end{thm}

\subsection{Proof:}The calculation of the contributions of components
$C$ such that
\[
q_S(S) \cap \Sg = \emptyset,
\]
is just as in Corollary~\ref{cor:Xsmooth}.  So by \ref{ass:Y}, we need
only treat the components of $C$ for which all curves of $S=S(C)$ lie in
a surface $Y_S\subset X_0$.   For such a $C$, and $x^i_S\in
Y_S\cap\Sg$, let  $Z\hra X_0$ be a fiber of $S$ over a point of
$D_{S'}^i$;  that is, 
$Z$ is a curve parameterized by $S'$ which passes through the node
$x^i_S$ of~$X_0$.  Since the embedding dimension of our threefold
singularity is four, we can choose local generators
\[
\{a_1, a_2\}\cup \{b_1.\dots,b_{\dim P -4}\}
\]
for the ideal $\I{Y_S}{P}$ near $x^i_S$ in such a way that the $b_k$
locally generate the ideal of a smooth fourfold containing~$X_0$.

By Lemma~\ref{lem:1.8}, there is an exact sequence
\begin{equation}\label{eq:IS}
0 \lra q^*\shd{\E}\tensor \O_{S} \xra{q^*s_0} \j\big/\j^2\tensor
 \O_{S} \lra \i_S/\i_S^2 \lra 0
\end{equation}
whence the exact sequence
\begin{multline}\label{eq:4term}
0 \lra \hom_{S}\left(\i_S\big/\i_S^2,\O_{S}\right) \lra
 \hom_{S}\left(\j_S\big/\j_S^2,\O_{S}\right) \lra \\
 \lra q_S^*\E \lra \ext_{S}^1\left(\i_S\big/\i_S^2,\O_{S}\right) \lra 0.
\end{multline}

Now since $X_0$ is the zero-scheme of~$s_0$, there is some
trivialization of $\E$ near 
$x^i_S$ with respect to which the local expression for $s_0$
is
\[
s_0 = \left(a_1f_1+a_2f_2, b_1,\dots,b_{\dim P -4}\right),
\]
with $\{f_1,f_2\}$ cutting out the point $x^i_S$ in~$Y_S$.
So under the map
\[
\hom_Z\left(\I{Z}{P}\big/\I{Z}{P}^2\tensor\O_Z,\O_Z\right)
\lra \E\tensor \O_Z,
\]
the images of the $((\dim P)-4)$ homomorphisms
\[
b_j\mapsto
\begin{cases}
1 \quad &\text{if $j=k$} \\
0 \quad &\text{otherwise}
\end{cases}
\]
generate a codimension one subspace of the vector
space~$\rest{E}{x^i_S}$.  We denote the one-dimensional quotient as
$\Co_i$.

Now the mapping
\[
\hom_{S}\left(\j_S\big/\j_S^2,\O_{S}\right) \lra q_S^*\E 
\]
is surjective away from $q_S^{-1}(\Sg)$, and since the generators
$a_1$ and $a_2$ must give sections of $\E\tensor\O_Z$ which
vanish at~$x^i_S$, we conclude from 
\eqref{eq:4term} that there is a natural isomorphism
\begin{equation}\label{eq:ext}
\bigoplus_i q_S^* \Co_i \isoa
 \ext_{S}^1\left(\i_S\big/\i_S^2,\O_{S}\right). 
\end{equation}

Recall that we have a functor
\[
\fn{T}^S := p^S_*\circ \hom_{S}\left(\blank,\O_{S}\right).
\]
Applying $\R\fn{T}^S$ to the sequence \eqref{eq:IS}. we obtain an exact
sequence
\begin{equation}
0 \lra \fn{T}^S\left(\i_S\big/\i_S^2\right) \lra \T_{J'}\tensor\O_{S'}
 \lra \V\tensor\O_{S'} \lra \R^1\fn{T}^S\left(\i_S\big/\i_S^2\right) \lra
 0.
\end{equation}
Because the fiber dimension of $p^S$ is one, the Grothendieck spectral
sequence for $\fn{T}^S$ degenerates at $E_2$.  We therefore obtain an
exact sequence:
\begin{multline}\label{eq:RT2}
0 \lra
\R^1p^S_*\left(\hom_{S}\left(\i_S\big/\i_S^2,\O_{S}\right)\right) \lra\\
 \lra \R^1\fn{T}^S\left(\i_S\big/\i_S^2\right) \lra
 p^S_*\ext^1_{S}\left(\i_S\big/\i_S^2,\O_{S}\right) \lra 0.
\end{multline}
Thus, we have the commutative diagram of $\O_{S'}$-modules
\begin{equation}\label{eq:diagram}
\begin{CD}
0 @. 0 @. {}\\
@VVV    @VVV @.\\
\sh{W}_{S'} 
 @>>> \R^1p^S_*\left(\hom_{S}\left(\i_S\big/\i_S^2,\O_{S}\right)\right)
 @>>> 0\\
@VVV @VVV @.\\
\V\tensor\O_{S'} @>>> \R^1\fn{T}^S\left(\i_S\big/\i_S^2\right) @>>> 0\\
@VV{\epsilon}V @VVV @.\\
\bigoplus_i\O_{D_{S'}^i} @>{\overset{\text{\eqref{eq:ext}}}{\hiso}}>>
 p^S_*\ext^1_{S}\left(\i_S\big/\i_S^2,\O_{S}\right) @. {}\\
@VVV @VVV @.\\
0 @. 0 @. {}
\end{CD}
\end{equation}
with exact rows and columns (where $\sh{W}_{S'}:=\ker(\epsilon)$).

Let $y'\in S'$;  we have a map
\[
\epsilon_{y'}\colon V_{y'}:= (\V\big/\m_{y',S'}\V) \to \bigoplus_i
 \left(\O_{D^i_{S'}}\big/\m_{y',S'}\O_{D^i_{S'}}\right) 
\]
induced from $\epsilon$ in \eqref{eq:diagram}.  Referring to
\eqref{eq:cones}, the condition \eqref{eq:Obs} says that, under the
geometric inclusion $C\subset \rest{V}{S'}$ we have:
\begin{equation}\label{eq:Ckernel}
\epsilon_{y'}(C \cap V_{y'}) =0.
\end{equation}
(It is here that we use the fact that the rays of the normal cone
correspond to curvilinear obstructions just as in the proof of
Lemma~\ref{lem:Cp}.) 

We next claim that there is a natural surjection
\begin{equation}\label{eq:surj}
\R^1 p^S_*\left(\hom_{S}\left(\i_S\big/\i_S^2,\O_{S}\right)\right)
 \lra \Omega^1_{S'}\tensor \K_{S'}^{-1}.
\end{equation}
To see this, let
\begin{equation}\label{eq:U}
U:=\left(S\setminus q_S^{-1}(\Sg)\right)\xra{j} S
\end{equation}
be the natural inclusion.  Then referring to
\ref{sec:adjunction}, we have morphisms 
\begin{equation}\label{eq:j}
\begin{split}
\omega_{S/S'}\tensor \left(\i_S\big/\i_S^2\right) &\lra
j_*j^*\left(\omega_{S/S'}\tensor \left(\i_S\big/\i_S^2\right)\right)\\
&\lra j_*j^*\left(\hom_{S}\left(\i_S\big/\i_S^2,\omega_{S/S'}\tensor
 \Lambda^2\left(\i_S\big/\i_S^2\right)\right)\right)\\
&\lra\hom_{S}\left(\i_S\big/\i_S^2,j_*j^*q_S^*\omega_{X_0}\right)\\
&\longleftarrow \hom_{S}\left(\i_S\big/\i_S^2,q_S^*\omega_{X_0}\right)\\
&\overset{\hiso}{\longleftarrow}
 \hom_{S}\left(\i_S\big/\i_S^2,\O_{S}\right)\tensor (p^S)^*\K_{S'}
\end{split}
\end{equation}
all of which restrict to isomorphisms over~$U$.  Moreover, the fiber
dimension of $p^S$ is one, so we have surjections 
\[
\minCDarrowwidth1pc
\begin{CD}
\R^1p^S_*\left(\omega_{S/S'}\tensor
 \left(\i_S\big/\i_S^2\right)\right) @. {} \\ 
@VVV @.\\
\R^1p^S_*\left(\hom_{S}\left(\i_S\big/\i_S^2,
 j_*j^*q_S^*\omega_{X_0}\right)\right) @<<< 
 \R^1p^S_*\left(\hom_{S}\left(\i_S\big/\i_S^2,\O_{S}\right)\right)
 \tensor\K_{S'} 
\end{CD}
\]
which are isomorphisms modulo torsion supported on~$D_{S'}$.  But by
Lemma~\ref{lem:1.9}, there is a surjection
\[
\R^1p^S_*\left(\omega_{S/S'}\tensor
 \left(\i_S\big/\i_S^2\right)\right) \isoa \Omega_{I'}^1\tensor\O_{S'}
 \xra{\mathrm{rest.}} \Omega^1_{S'}
\]
and thus, since $\Omega_{S'}^1$ is locally free, we can deduce the
existence of the surjection~\eqref{eq:surj}, and therefore, from
\eqref{eq:diagram}, as surjection
\[
\sh{W}_{S'} \to \Omega_{S'}^1\tensor\K_{S'}^{-1}.
\]   
Set
\[
\Q_{Q'} := 
\]

Referring to \eqref{eq:diagram}, let
\[
\Q_{S'} := \frac{\V\tensor\O_{S'}}{\ker\left(\sh{W}_{S'} \to
\Omega^1_{S'}\tensor\K_{S'}^{-1}\right)}. 
\]
Now a diagram chase completes \eqref{eq:diagram} to a commutative
diagram 
\begin{equation}\label{eq:complete}
\begin{CD}
0 @. 0 @. 0\\
@VVV  @VVV @VVV\\
\sh{W}_{S'} @>>> 
 \R^1p^S_*\left(\hom_{S}\left(\i_S\big/\i_S^2,\O_{S}\right)\right)
 @>{\eqref{eq:surj}}>>  \Omega_{S'}^1\tensor\K_{S'}^{-1}\\
@VVV  @VVV @VVV\\
\V\tensor\O_{S'} @>>> \R^1\fn{T}^S\left(\i_S\big/\i_S^2\right) @>>> \Q_{S'}\\
@VV{\epsilon}V @VVV @VVV\\
\bigoplus_i\O_{D_{S'}^i} @>{\overset{\text{\eqref{eq:ext}}}{\hiso}}>>
 p^S_*\ext^1_{S}\left(\i_S\big/\i_S^2,\O_{S}\right) @>{\hiso}>>
 \bigoplus_i\O_{D_{S'}^i}\\
@VVV @VVV @VVV\\
0 @. 0 @. 0
\end{CD}
\end{equation}
with exact columns and horizontal epimorphisms.  

Finally consider \ref{ass:nc}.  In case~1, $\sh{W}_{S'}$ is locally free
because the normal-crossing divisor $D_{S'}$ imposes independent
conditions:  if $y'\in S'$ lies in the divisors
$D_{S'}^{i_1},\dots,D_{S'}^{i_r}$, then near $y'$, the left-hand
column of \eqref{eq:complete} looks like
\[
0 \lra \bigoplus_{j=1}^{r} \O_{S'}\left(-D_{S'}^{i_j}\right) \oplus
\O_{S'}^{r'} \lra \O_{S'}^{r+r'} \lra
\bigoplus_{j=1}^{r}\O_{D_{S'}^{i_j}} \lra 0.
\]

Then away from $D_{S'}$,
the cone component $C$ must be, as in the smooth case, exactly the
geometric kernel of the epimorphism of vector bundles
\[
W_{S'} \lra \shd{T}_{S'}\tensor K_{S'}^{-1}.
\]
By \eqref{eq:Ckernel} therefore, $C$ must be a geometric sub-bundle of
$\rest{V}{S'}$, so that $\Q_{S'}$ is locally free.  Then, proceeding just
as in the proof of Corollary~\ref{cor:Xsmooth}, we have that
\[
\gamma(I') = \sum_{C\in \Ccomp}m(C) \ctop(\Q_{S'})\cap [S'(C)].
\]

In case~2 of \ref{ass:nc}, i.e., if $S'=I'$ is smooth,
$\R^1p^S_*\left(\hom_{S}\left(\i_S\big/\i_S^2,\O_{S}\right)\right)$ is
locally free by semicontinuity and the local theory of the Hilbert
scheme.  Consequently, since the surjection of \eqref{eq:surj} is
generically injective (Proposition~\ref{prop:OmegaI'}), it is an
isomorphism.  Thus,
\[
\R^1\fn{T}^S\left(\i_S\big/\i_S^2\right) \xra{\hiso} \Q_{S'}.
\]

Now $C$ is just the geometric normal bundle of $S'=I'$ in $J'$ and
keeping in mind the infinitesimal properties of the Hilbert scheme,
\eqref{eq:RT} (with $A'=S'=I'$) realizes
$\R^1\fn{T}^S\left(\i_S\big/\i_S^2\right)$ as the excess normal bundle
of the diagram 
\[
\begin{CD}
{} @. J'\\
@. @VV{\sigma_0}V\\
J' @>{z_V}>> V
\end{CD}
\]
The desired formula is now an immediate consequence of the excess
intersection formula.\qed

\subsection*{}Finally, we have the following analogue of
Proposition~\ref{prop:OmegaI'}: 
\begin{thm}\label{thm:log}
If, in Theorem~\ref{thm:nodes}, $I'$ itself is smooth and, for the
generic point $\p^i$ of each component $D^i$ of~$D=D_{I'}$, the curves
$I_{\p^i}$ are locally smooth at the node $x^i$, then
\[
\R^1\fn{T}(\i\big/\i^2) \iso \K_{I'}^{-1}\tensor\Omega^1_{I'}[\log D]
= \Q_{I'}. 
\]
\end{thm}

\subsection{Proof:}We refer to the proof of Theorem~\ref{thm:nodes}
in the case $S'=I'$ (in which case $\i_S=\i$ and $p^S=p^0$).  We have
seen that the kernel of the surjection 
\begin{equation}\label{eq:surj2}
\R^1p^0_*\left(\hom_I\left(\i\big/\i^2,\O_I\right)\right) \to
\Omega^1_{I'}\tensor\K_{I'}^{-1} 
\end{equation}
given in \eqref{eq:surj} is supported along~$D$.  On the other hand, by
\eqref{eq:RT2}, the domain of \eqref{eq:surj2} is a subsheaf of the
locally free sheaf $\Q_{I'}$ and hence, so is its kernel.  Since $I'$
is smooth, it follows that \eqref{eq:surj2} 
must in fact be an isomorphism. 
So, by \eqref{eq:complete},
\[
\R^1\fn{T}(\i\big/\i^2) = \Q_{I'}
\]
and is therefore locally free.  Now apply the functor
$\R\hom_{I'}(\blank,\O_{I'})$ to the exact sequence \eqref{eq:RT2} to
obtain the exact sequence
\begin{equation}\label{eq:log}
0\lra \shd{\Q}_{I'} \lra \T_{I'}\tensor\K_{I'} \xra{\tau}
\bigoplus_i\O_{D^i}(D^i) \lra 0 
\end{equation}
since, by the standard divisorial exact sequence,
\[
\ext^1_{I'}(\O_{D^i},\O_{I'}) \iso \O_{D^i}(D^i).
\]
Thus it suffices to show that, at a general point $y'$ of any
fixed~$D^i$, the sections of $\ker\tau$ lie in
$\K_{I'}\tensor\T_{I'}[\log D^i]$. 
Choose coordinates $y'_0,\dots$ on $I'$ near $y'=0$ such that 
\[ 
D^i = \{y'_0=0\}.
\] 

We refer to the \eqref{eq:U} and \eqref{eq:j} in the proof of
Theorem~\ref{thm:nodes}. We analyze the local behavior of the mapping
\begin{equation}\label{eq:mapping}
\left(\i\big/\i^2\right)\tensor\omega_{I/I'} \lra
\hom_I(\i\big/\i^2,j_*j^*q_0^*\omega_{X_0}) 
\end{equation}
in terms of analytic local coordinates
\[
\{c_1,c_2\}\cup\{a_1,a_2\}\cup\{b_1,\dots,b_{\dim P -4}\}
\]
defined on a small analytic neighborhood of the node~$x^i\in P$.
(Compare \eqref{eq:rel_adj}).  Since $j$ is the inclusion of an
open subset with complement of codimension~$2$, then locally
\[
j_*j^*q^*\omega_{X_0}\iso q^*\omega_{X_0}.
\]
We may assume that upon restricting to the smooth analytic fourfold
$F$ given by 
\[
\{b_1=\dots=b_{\dim P -4}=0\},
\]
$X_0$ has local equation 
\[
a_1c_1+a_2c_2 = 0,
\]
the curve $I_0$ is locally given by 
\[
c_1=a_1=a_2=0,
\]
and the incidence scheme $I$ is given locally in $I'\times F$ by
\begin{align*}
a_1=a_2&=0\\
c_1+c_1f_0(c_1,c_2)+c_2g_0(c_1,c_2) &= ky'_0,
\end{align*}
with $k\ne0$ and $f_0,g_0$ vanishing at~$(0,0)$.  Also, $\omega_{X_0}$
has local generator
\begin{multline}\label{eq:Res}
\Res_{X_0}\frac{da_1\wedge da_2\wedge dc_1\wedge dc_2}{a_1c_1+a_2c_2} =\\
 \pm\frac{da_1\wedge dc_1\wedge dc_2}{c_2}
= \pm\frac{da_1\wedge da_2\wedge dc_2}{a_1}
= \pm\frac{da_2\wedge dc_1\wedge dc_2}{c_1},
\end{multline}
whereas $\omega_{I/I'}$ has local generator~$dc_2$.
Thus, keeping \eqref{eq:Ktrivial} in mind, \eqref{eq:mapping}
specializes at the curve $I_0$ to the map 
\begin{equation}
\left(\I{I_0}{X_0}\big/\I{I_0}{X_0}^2\right)\tensor\omega_{X_0} \lra
 \hom_{I_0}\left(\I{I_0}{X_0}\big/\I{I_0}{X_0}^2,\O_{I_0}\right)
\end{equation}
given by
\begin{equation}\label{eq:special_mapping}
\begin{split}
a_i & \mapsto \begin{pmatrix}a_j & \mapsto &\pm da_i\wedge da_j\wedge
dc_2 \\
c_1 & \mapsto & \pm da_i\wedge dc_1\wedge dc_2\end{pmatrix}\\
c_1 & \mapsto \begin{pmatrix}a_i & \mapsto & \pm
da_i\wedge dc_1\wedge dc_2\\c_1 & \mapsto &0\end{pmatrix}.\\
\end{split}\end{equation}
But, by \eqref{eq:Res}, the image of \eqref{eq:special_mapping}
consist entirely of homomorphisms which vanish at the point~$x^i$.
Thus, near $y'=0$, the 
mapping 
\[
\R^1p^0_*(\omega_{I/I'}\tensor(\i\big/\i^2)) \lra
\R^1p^0_*(\hom_I(\i\big/\i^2,q_0^*\omega_{X_0}))
\]
is a surjection with torsion kernel, so that, from \eqref{eq:RT2} and
the projection formula we have the exact sequence
\[
0\lra\frac{\R^1p^0_*(\omega_{I/I'}\tensor(\i\big/\i^2))}{(\text{some
torsion sheaf})}\tensor\K_{I'}^{-1} \lra
\R^1\fn{T}(\i\big/\i^2)\tensor \lra p^0_*\ext^1_{I}(\i\big/\i^2,\O_I)  \lra 0.
\]
So, as in \eqref{eq:log}, we have an exact sequence
\[
0\to\shd{\left(\R^1\fn{T}\left(\i\big/\i^2\right)\right)}  \to
\shd{\left(\frac{\R^1p^0_*(\omega_{I/I'}\tensor(\i\big/\i^2))}
{(\text{some torsion sheaf})}\tensor\K_{I'}^{-1}\right)}\xra{\tau}
\bigoplus_i \O_{D^i}(D^i)\to 0,   
\]  
where the middle term is isomorphic to $\T_{I'}\tensor\K_{I'}$.
But by Verdier duality:
\begin{equation}\label{eq:Verdier3}
\begin{split}
\R^1\fn{T}(\i\big/\i^2) &=
 \R^1(p^0_*\circ\hom_I(\blank,\omega_{I/I'}))
 (\i\big/\i^2\tensor\omega_{I/I'}) \\ 
&= \hom_{I'}(p^0_*((\i\big/\i^2)\tensor\omega_{I/I'}),\O_{I'}).
\end{split}
\end{equation}
So, by local freeness and by \eqref{eq:log}, we obtain an exact sequence
\[
p^0_*((\i\big/\i^2)\tensor\omega_{I/I'}) \xra{\rho} \T_{I'}\tensor\K_{I'}
 \xra{\tau} \bigoplus_i\O_{D^i}(D^i) \lra 0,
\]
where the image of $\rho$ consists entirely of homomorphisms which
vanish at $x^i$ and so lie in~$\K_{I'}\tensor\T_{I'}[\log D^i]$.\qed 
\section{Applications}\label{sec:exas}
We conclude this paper with three applications of the formulas of
Corollary~\ref{cor:Xsmooth} and Theorem~\ref{thm:nodes}.

\subsection{Example:}  This application was suggested by A.~Bertram and
M.~Thaddeus.  Let $C$ be a hyperelliptic curve of genus~$4$ and let
\[
X_0 = C^{(3)}
\]
be the third symmetric power of ~$C$.  Embed $X_0$ in 
\[
P = C^{(7)}
\]
which is a \P[3] bundle over the Jacobian~$J(C)$.  Thus, rational
curves in $P$ are unobstructed and $X_0$ is the zero-scheme of a
section of
\[
\E = \L^{\oplus4},
\]
where $\L$ is the line bundle given by the divisor
\[
(\mathrm{basept.}+ C^{(6)}).
\]
Let 
\[
I' = W_3^1 \iso C
\]
be the Hilbert scheme of $g_3^1$'s on $C$ so that 
\[
I\subseteq I'\times C^{(3)}
\]
becomes the tautological \P[1]-fibration over~$I'$.  Under the
Abel-Jacobi map
\[
X_0 \to \Pic^3(C),
\]
the fibers of $I/I'$ are contracted to double points of the theta
divisor $\Theta$, which is itself the image of $X_0$.  These are
canonical singularities, so $\omega_{X_0}$ is the pullback
of~$\O(\Theta)$.  Thus
\[
(p^0)^*\left(\O(\Theta)\tensor \O_C\right) = q_0^*\omega_{X_0}
\]
so that Corollary~\ref{cor:Xsmooth} yields the well-known fact that
the number of $g_3^1$'s on a generic curve of genus four is
\[
\gamma(I') = c_1(\omega_C) - (\Theta\cdot C) = 6-4 =2.
\]

\subsection{Example:}  Our second application is again not new, being
the subject of~\cite{AK91}.  Let $X_0$ be the Fermat quintic
hypersurface in~$P=\P[4]$.  Let $I'$ be the Hilbert scheme of lines
in~$X_0$.  Following \cite{AK91}, $I'\red$ is the union of $50$ Fermat
quintic plane curves meeting transversely in pairs at $375$ points,
these points being exactly the flex point of the Fermat plane curves.
Using Pl\"ucker coordinates in the Grassmannian of~\P[4], one computes
that the local analytic structure of $I'$ away from the crossings is
given by
\[
\Co[x,y]\big/(y^2)
\]
while at each of the 375 crossings it is given by
\[
\Co[x,y]\big/(x^3y^2,x^2y^3).
\]
The components of the normal cone are computed from this local
analytic structure using the local primary decomposition
\[
(x^3y^2,x^2y^3) = (y^2)\cap(x^2)\cap(x,y)^5.
\]
One computes that the normal cone has one component of multiplicity $2$
over each Fermat quintic curve $F$ and one component of multiplicity
$5$ over each crossing point.  Hence, by Corollary~\ref{cor:Xsmooth}, the
number of lines on the general quintic threefold is
\[
\gamma(I') = 50\cdot 2\cdot c_1(\omega_F) +
 375\cdot5\cdot c_0(\omega_{\mathrm{pt.}}) = 2875.
\]

\subsection{Example}\label{ex:ciCY}  The setting for the final
 application was first considered \cite{Clemens83}, and later
in~\cite{Kley98}, \cite{EJS97} and \cite{Kley98a}. To the 
authors' best knowledge, however, the calculations below are new. 

Consider either of the following cases:  Choose
\begin{equation}\label{eq:41}
g_1,\alpha_2\in\Gamma(\P[4],\O(4))\quad\text{and}\quad
g_2,\alpha_1\in\Gamma(\P[4],\O(1)) 
\end{equation}
or
\begin{equation}\label{eq:32}
g_1,\alpha_2\in\Gamma(\P[4],\O(3))\quad\text{and}\quad
g_2,\alpha_1\in\Gamma(\P[4],\O(2)) 
\end{equation}
sufficiently general so that both the K3 surface
\[
Y:=\{g_1=g_2=0\}
\]
and the del Pezzo surface
\[
S:=\{g_2=\alpha_1=0\}
\]
are smooth, and such that the quintic threefold
\[
X_0:=\{\alpha_1g_1 + \alpha_2g_2 = 0\}
\]
has only ordinary nodes, all of which---$16$ in case \eqref{eq:41} and
$36$ in case \eqref{eq:32}---lie on~$S\cap Y$.  In \cite{Kley98}, it is
shown that, despite the existence of the nodes,
\[
\N{Y}{X_0} = \omega_{Y}\quad\text{and}\quad\N{S}{X_0} = \omega_S.
\] 
For any curve $C$ in $S$, we have the exact sequence
\[
0\lra \H^0(C,\N{C}{S})\lra\H^0(C,\N{C}{X_0})\lra\H^0(C,\N{S}{X_0}\tensor\O_C),
\]
and, since $\shd{\omega_S}$ is ample and $\h^1(\O_S) = 0$,
\[
\h^0(\N{C}{X_0}) = \h^0(\N{C}{S})  = \h^0(\O_S(C))-1.
\]
Thus the linear system
\[
I':=\abs{\O_S(C)}
\]
is a connected component of the Hilbert scheme of~$X_0$.  If the $g_i$
and the $\alpha_i$ are sufficiently general, the divisor of curves
passing through at least one node is a simple normal-crossing divisor
consisting of hyperplanes.  So we may apply Theorem~\ref{thm:log} as
long as 
\[
\H^1(C,\N{C}{\P[4]}) = 0
\]
for all $C$ in~$I'$.

For example, in case \eqref{eq:41}, the lines in the plane $S$
contribute
\[
\gamma(I') = c_2(\Omega^1_{\P[2]}[\log 16\P[1]])\cap[\P[2]] = 91
\]
lines to the general quintic threefold~$X$, the conics in~$S$ contribute
\[
c_5(\Omega^1_{\P[5]}[\log 16\P[4]])\cap[\P[5]] = 2002
\]
conics to~$X$, and the cubic curves in $S$ contribute
\[
c_9(\Omega^1_{\P[9]}[\log 16\P[]])\cap[\P[9]] = 2002
\]
cubic elliptic curves to~$X$.

In case \eqref{eq:32},
\[
S\iso \Bl_{\text{$5$ pts.}}\P[2]
\]
and the lines in \P[2] contribute
\[
c_2(\Omega^1_{\P[2]}[\log 36\P[1]])\cap[\P[2]] = 595.
\]
twisted cubics to~$X$, and the hyperplane sections of $S$ contribute
\[
c_4(\Omega^1_{\P[4]}[\log 36\P[3]])\cap[\P[4]] = 46,\!376
\]
degree~$4$ elliptic curves to~$X$.

\medskip
Note that one can perform analogous constructions and computation on the other
types of complete-intersection Calabi-Yau threefolds.  Moreover, in
any of these cases, the contributions of curves lying on the 
K3 surface $Y$ can be calculated;  see
\cite{Kley98} (curves of genus~$1$) and \cite{Kley98a} (where the
results of the present paper are applied to curves of higher genus).
In fact, it is shown there that if $Y$ has Picard number~$2$, then the
curves coming from the primitive linear system on $Y$ not generated by
the hyperplane sections contribute only geometrically rigid curves to
a general deformation of~$X_0$.

\bibliographystyle{amsplain}
\bibliography{mybib}


\end{document}